\documentstyle{amsppt}
\baselineskip18pt
\magnification=\magstep1
\pagewidth{30pc}
\pageheight{45pc}

\hyphenation{co-deter-min-ant co-deter-min-ants pa-ra-met-rised
pre-print pro-pa-gat-ing pro-pa-gate
fel-low-ship Cox-et-er dis-trib-ut-ive}
\def\leaderfill{\leaders\hbox to 1em{\hss.\hss}\hfill}

\def\idest{i.e.,\ }

\def\a{{\alpha}}
\def\be{{\beta}}
\def\g{{\gamma}}

\def\bn{{\bold n}}

\def\bs{{\bold s}}

\def\bx{{\bold x}}

\def\b0{\text{\bf 0}}

\def\ra{{\ \longrightarrow \ }}

\def\rank{\text{\rm \, rank}}

\def\real{{\Bbb R}}
\def\complex{{\Bbb C}}
\def\zed{{\Bbb Z}}
\def\kyu{{\Bbb Q}}

\def\Im{\text{\rm Im}}

\def\pd{\partial}
\def\boxit#1{\vbox{\hrule\hbox{\vrule \kern3pt
\vbox{\kern3pt\hbox{#1}\kern3pt}\kern3pt\vrule}\hrule}}
\def\rabbit{\vbox{\hbox{\kern0pt
\vbox{\kern0pt{\hbox{---}}\kern3.5pt}}}}

\def\tableau#1{
        \hbox {
                \hskip -10pt plus0pt minus0pt
                \raise\baselineskip\hbox{
                \offinterlineskip
                \hbox{#1}}
                \hskip0.25em
        }
}

\def\tabCol#1{
\hbox{\vtop{\hrule
\halign{\strut\vrule\hskip0.5em##\hskip0.5em\hfill\vrule\cr\lower0pt
\hbox\bgroup$#1$\egroup \cr}
\hrule
} } \hskip -10.5pt plus0pt minus0pt}

\def\CR{
        $\egroup\cr
        \noalign{\hrule}
        \lower0pt\hbox\bgroup$
}



\def\blank#1#2{
\hbox to #1{\hfill \vbox to #2{\vfill}}
}


\def\strut{\vrule height10pt depth5pt width0pt}

\def\rank{\text{\rm rank}}

\def\hcn{{\Gamma_n}}
\def\oldpsi{{V}}

\def\ven{{V^e_n}}
\def\von{{V^o_n}}
\def\oldbv{{v}}
\def\boz{{\underline{0}}}
\def\boo{{\underline{1}}}
\def\boy{{\bold y}}

\def\secz{1}
\def\secaz{2}
\def\secba{3}
\def\secc{4}
\def\secd{5}
\def\sece{6}

\topmatter
\title Morse matchings on polytopes \endtitle

\author R.M. Green and Jacob T. Harper \endauthor
\affil Department of Mathematics \\ University of Colorado \\
Campus Box 395 \\ Boulder, CO  80309-0395 \\ USA \\ {\it  E-mail:}
rmg\@euclid.colorado.edu,\ Jacob.Harper\@colorado.edu \\
\newline
\endaffil


\abstract
We show how to construct homology bases for certain CW complexes in terms
of discrete Morse theory and cellular homology.  We apply this technique to
study certain subcomplexes of the half cube polytope studied in previous 
works.  This involves constructing explicit complete acyclic Morse matchings 
on the face lattice of the half cube; this procedure may be of independent
interest for other highly symmetric polytopes.
\endabstract

\endtopmatter

\centerline{\bf Preliminary version, draft 2}

\head \secz. Introduction \endhead

As a graph, the $n$-dimensional hypercube is bipartite and connected.
This induces a partition of its vertex set $\oldpsi = \oldpsi_n = \{\pm 1\}^n$
into two pieces, $\oldpsi^e \cup \oldpsi^o = \oldpsi_n^e \cup \oldpsi_n^o$, 
where $\oldpsi_n^e$ (respectively, $\oldpsi_n^o$) consists of those vertices 
whose coordinates contain an even (respectively, odd) number of occurrences
of $-1$.  We define the {\it half cube}, $\hcn$, to be the convex hull of the
$2^{n-1}$ points in $\oldpsi_n^e$.  Using $\oldpsi_n^o$ in place of 
$\oldpsi_n^e$ in this 
construction gives rise to an isometric copy of the half cube.

In a previous work \cite{{\bf 10}}, 
the first author classified the faces of the half cube and explained how they
assemble naturally into a regular CW complex, $C_n$, which is homeomorphic 
to a ball (see Theorem \secc.2).  Furthermore, for each 
$3 \leq k \leq n$, there is an interesting subcomplex $C_{n, k}$ of $C_n$
obtained by deleting the interiors of all the half cube shaped faces of 
dimensions $l \geq k$.  We also showed in \cite{{\bf 10}, Theorem 3.3.2} 
that the reduced homology of
$C_{n, k}$ is free over $\zed$ and concentrated in degree $k-1$.

The Coxeter group $W(D_n)$ acts naturally on the $(k-1)$-st homology of
$C_{n, k}$, and we computed the character of this representation 
(over $\complex$) in \cite{{\bf 11}, Theorem 4.4}.  The group $W(D_n)$ has
two parabolic subgroups that are isomorphic to the symmetric group
$S_n$, so the homology representations become $S_n$-modules by restriction.
We showed in \cite{{\bf 11}, Theorem 4.7} that the resulting representations 
of $S_n$ are equivalent to the representation of $S_n$ on the $(k-2)$-nd 
homology of the complement of the $k$-equal real hyperplane arrangement.

If $k = n-1$, the complex $C_{n, k}$ is the boundary complex of the half
cube, and is therefore shellable by a well-known result of Bruggesser--Mani
\cite{{\bf 4}}.  If $k < n-1$, then the fact that $C_{n, k}$ has nonzero
$(k-1)$-st homology is an obstruction to shellability, which means that
we cannot use the machinery of shellability to produce a homology basis
for $C_{n, k}$.

The first main result of this paper (Theorem \secba.7) shows how to
use cellular homology and discrete Morse theory to construct an explicit
integral homology basis for certain kinds of regular CW complexes, of
which the complexes $C_{n, k}$ are motivating examples.  In order to apply
the theorem to a CW complex $Y$, one starts with a complete acyclic Morse 
matching $V$ of a CW complex $X$ that contains $Y$ as a subcomplex.  Under
certain mild additional hypotheses, which are satisfied by $C_{n, k}$,
the theorem produces an explicit set of boundaries in $X$ that induce a
homology basis for $Y$.

In \S\secc, we construct a complete acyclic matching on the face lattice of
the half cube.  It is not a surprise that such a matching exists, as this
follows abstractly from the shellability of the boundary complex of a 
polytope.  Our motivation is to construct an explicit matching that works
for half cubes of arbitrary dimension $n \geq 4$.  We prove in Theorem \secd.8
that this is a complete acyclic matching on the faces of the half cube
(together with the empty face).

Let $b_{n, k}$ be the $(k-1)$-st Betti number of the subcomplex $C_{n, k}$.
The numbers $b_{n, k}$ appear as sequence A119258 in Sloane's online
encyclopedia \cite{{\bf 16}}.  Various explicit expressions for these
numbers are given in \cite{{\bf 15}}, some of which are sums of products of 
positive
integers.  One of these appears in work of Bj\"orner--Welker \cite{{\bf 3}}, 
who study the numbers $b_{n, k}$ in the context of hyperplane arrangements.
They prove that $$
b_{n, k}  = \sum_{i = k}^n {n \choose i} {{i-1} \choose {k-1}}. 
\eqno{(\secz.1)}
$$  There is a representation theoretic explanation for this: the terms in
the sum correspond to the dimensions of the irreducible constituents of
the representation of the Coxeter group $W(D_n)$ on the reduced homology
of $C_{n, k}$; this was shown in \cite{{\bf 11}}.

Another formula for $b_{n, k}$ is $$
b_{n,k} = \sum_{i = 1}^n 2^{i-k} {{i-1} \choose {k-1}}; \eqno{(\secz.2)}
$$ this is a straightforward generalization of a result of Barcelo and
Smith \cite{{\bf 1}}, who study the case $k = 3$ in the context of combinatorial
homotopy theory (``$A$-theory'').  In Theorem \sece.2, we construct an 
explicit homology basis for $C_{n, k}$ in terms of cellular homology; this 
basis has the property that when it is enumerated in the obvious way, we 
recover equation (\secz.2).

We believe that the quest for explicit Morse matchings on the face lattices
of polytopes is an aesthetically pleasing goal in its own right, similar to
the discovery of an explicit shelling order on the faces of a polytope.  
It would be interesting to find such matchings for other polytopes.  For 
some, such as the hypercube and the simplex, this is a fairly easy exercise.
Others, such as the hypersimplex, present about the same level of
difficulty as the half cube; this is described in the second author's thesis
\cite{{\bf 12}} and will be published separately.  It would be very interesting 
to have such a description for the permutahedron, some of whose subcomplexes 
are known to have important topological properties \cite{{\bf 2}, Theorem 2.4}.


\head \secaz. Discrete Morse theory for CW complexes \endhead

We first recall the definition of a finite regular CW complex, following
\cite{{\bf 14}, \S8}.

If $X$ and $Y$ are topological spaces with $A \subset X$ and $B \subset Y$,
we define a continuous map $g : (X, A) \ra (Y, B)$ to be a continuous
map $g : X \ra Y$ such that $g(A) \subseteq B$.  If, furthermore,
$g|_{X - A} : X - A \ra Y - B$ is a homeomorphism, we call $G$ a {\it relative
homeomorphism}.

An {\it $n$-cell}, $e = e^n$ is a homeomorphic copy of the open $n$-disk
$D^n - S^{n-1}$, where $D^n$ is the closed unit ball in Euclidean $n$-space
and $S^{n-1}$ is its boundary, the unit $(n-1)$-sphere.  We call $e$ a 
{\it cell} if it is an $n$-cell for some $n$.

If a topological space $X$ is a disjoint union of cells $X = \bigcup \{e :
e \in E\}$, then for each $k \geq 0$, we define the $k$-skeleton
$X^{(k)}$ of $X$ by $$
X^{(k)} = \bigcup \{e \in E : \dim(e) \leq k\}
.$$

The CW complexes we consider in this paper are all finite, which means that
we can give the following abbreviated definition.

\definition{Definition \secaz.1}
A {\it CW complex} is an ordered triple $(X, E, \Phi)$, where $X$ is a 
Hausdorff space, $E$ is a family of cells in $X$, and $\{\Phi_e : e \in E\}$
is a family of maps, such that 
\item{(i)}{$X = \bigcup \{e : e \in E\}$ is a disjoint union;}
\item{(ii)}{for each $k$-cell $e \in E$, the map $\Phi_e : 
(D^k, S^{k-1}) \ra (e \cup X^{(k-1)}, X^{(k-1)})$ is a relative homeomorphism.}

A {\it subcomplex} of the CW complex $(X, E, \Phi)$ is a triple $(|E'|, E',
\Phi')$, where $E' \subset E$, $$
|E'| := \bigcup \{ e : e \in E' \} \subset X
,$$ $\Phi' = \{\Phi_e : e \in E'\}$, and $\Im\  \Phi_e \subset |E'|$ for every
$e \in E'$.

The complexes considered here have the 
property that the maps $\Phi_e$ (regarded as mapping to their images) are 
all homeomorphisms.  Such CW complexes are called {\it regular}.

An {\it oriented CW complex} is a CW complex together with a choice of
orientation for each cell.
\enddefinition

Cellular homology is a version of singular homology that is particularly
convenient in the context of regular CW complexes.  For our purposes, it
is convenient to describe cellular homology in terms of intersection numbers
as follows.
If $e_{\a}^n$ is an $n$-cell and $e_{\be}^{n-1}$ is an $(n-1)$-cell of the
same CW complex $X$, then the {\it incidence number} $[e_\a^n : e_\be^{n-1}]$
is defined in \cite{{\bf 9}, \S2.5} as the degree of a certain map.  It follows
that the incidence number is an integer.  A key property for our purposes
is the following.

\proclaim{Proposition \secaz.2}
If $X$ is an oriented regular CW complex then the intersection number
$[e_\a^n : e_\be^{n-1}]$ is equal to $\pm 1$ if $e_\be^{n-1}$ is a face
of $e_\a^n$, and is equal to $0$ otherwise.
\endproclaim

\demo{Proof}
This is \cite{{\bf 9}, Proposition 5.3.10}.
\qed\enddemo

To define the cellular homology of a CW complex $X$ over a ring $R$, we 
introduce, for each integer $n \geq 0$, the $n$-chains of $X$.
This is the free $R$-module with a basis indexed by all the $n$-cells,
$e_\a^n$; by abuse of notation, we will identify the basis elements with
the cells, having fixed once and for all on an orientation for each cell.  
The boundary map $\pd = \pd_n : C_n(X; R) \ra C_{n-1}(X; R)$ is then
defined to be the $R$-module homomorphism for which $$
\pd(e_\a^n) = \sum_\be [e_\a^n : e_\be^{n-1}] e_\be^{n-1}
.$$  It can be shown that $\pd \circ \pd = 0$.  The homology of the complex
$C_\bullet$ is the {\it cellular homology} of $X$ over $R$.  It is
convenient for some purposes to introduce a unique $-1$-cell $e_\a^{-1}$;
this gives rise to {\it reduced} cellular homology.

Discrete Morse theory, which was introduced by Forman \cite{{\bf 7}}, is a
combinatorial technique for computing the homology of CW complexes.  
By building on work of Chari \cite{{\bf 6}}, Forman later produced a version
of discrete Morse theory based on acyclic matchings in Hasse diagrams
\cite{{\bf 8}}.  This version of the theory plays a key role in computing
the homology of $C_{n, k}$.

\definition{Definition \secaz.3}
Let $K$ be a finite regular CW complex.  
A {\it discrete vector field} on $K$ is
a collection of pairs of cells $(K_1, K_2)$ such that 
\item{(i)}{$K_1$ is a face of $K_2$ of codimension $1$ and}
\item{(ii)}{every cell of $K$ lies in at most one such pair.}

We call a cell of $K$ {\it paired} if it lies in (a unique) one of the above
pairs, and {\it unpaired} otherwise.  If $(K_1, K_2)$ is a pair of the matching
as above, we say that $K_1$ is an {\it upward matching} face and that $K_2$
is a {\it downward matching} face.

If $V$ is a discrete vector field on a regular CW complex $K$, we define
a $V$-path to be a sequence of cells $$
\a_0, \be_0, \a_1, \be_1, \a_2, \ldots, \be_r, \a_{r+1}
$$ such that for each $i = 0, \ldots, r$, (a) each of $\a_i$ and $\a_{i+1}$
is a codimension $1$ face of $\be_i$, 
(b) each $(\a_i, \be_i)$ belongs to $V$ and 
(c) $\a_i \ne \a_{i+1}$ for all $0 \leq i \leq r$.  If
$r \geq 0$, we call the $V$-path {\it nontrivial}, and if $\a_0 = \a_{r+1}$,
we call the $V$-path {\it closed}.  Note that all the faces $\a_i$ have
the same dimension, $p$ say, and all the faces $\be_i$ have dimension $p+1$.
\enddefinition

Let $P$ be the set of cells of $K$, together with the empty cell
$\emptyset$, which we consider to be a cell of dimension $-1$.  Denote
the set of $k$-cells of $K$ by $P^k$.  The set $P$
becomes a partially ordered set under inclusion.  Let $H$ be the Hasse diagram
of this partial order.  We regard $H$ as a directed graph, in which all edges
point towards cells of larger dimension.

Suppose now that $V$ is a discrete vector field on $K$.  We define $H(V)$ to
be the directed graph obtained from $H$ by reversing the direction of an arrow
if and only if it joins two cells $K_1 \subset K_2$ for which $(K_1, K_2)$
is one of the pairs of $V$.  If the graph $H(V)$ has no directed cycles, we
call $V$ an {\it acyclic (partial) matching} of the Hasse diagram of $K$.

\proclaim{Theorem \secaz.4 (Forman)}
Let $V$ be a discrete vector field on a regular CW complex $K$.
\item{\rm (i)}{There are no nontrivial closed $V$-paths if and only if $V$ is
an acyclic matching of the Hasse diagram of $K$.}
\item{\rm (ii)}{Suppose that $V$ is an acyclic partial matching of the Hasse 
diagram of $K$ in which the empty set is unpaired.  Let $u_p$ denote the number
of unpaired $p$-cells.  Then $K$ is homotopic to a CW complex
with exactly $u_p$ cells of dimension $p$ for each $p \geq 0$.}
\endproclaim

\demo{Proof}
Part (i) is \cite{{\bf 8}, Theorem 6.2} and part (ii) is 
\cite{{\bf 8}, Theorem 6.3}.
\qed\enddemo



\head \secba. Homology bases for subcomplexes \endhead

\definition{Definition \secba.1}
Let $X$ be a finite regular oriented CW complex 
and let $V$ be an acyclic partial matching on $K$.  
Recall that $P^{i}$ is the set of $i$-cells  of $K$.
For each $k$, let $V_k = V \cap (P^k \times P^{k+1})$.
Define $$
d_{V, k} = \{K \in P^{k+1} : (K_1, K) \in V_k 
\text{\ for some\ } K_1 \in P^k\}
$$ and $$
e_{V, k} = \{K \in P^k : (K, K_2) \in V_k
\text{\ for some\ } K_2 \in P^{k+1}\}
.$$  Let $D_{V, k}(X; R)$ (respectively, $E_{V, k}(X; R)$) 
be the free $R$-module on $d_{V, k}$ (respectively, $e_{V, k}$).
Let $\pd_{V, k} : D_{V, k}(X; R) \ra E_{V, k}(X; R)$ be the $R$-module
homomorphism defined by $$
\pd_{V, k}(e_\a^{k+1}) = \sum_{\be \in e_{V, k}} 
[e_\a^{k+1} : e_\be^k] e_\be^k
$$ for each $e_\a^{k+1}$.

If $e, e' \in e_{V, k}$, we write $e' \prec e$ if both $(e, d) \in V_k$ 
and $e'$ lies in the boundary of $d$.  Let $\leq_{e, k}$ be the relation on
$e_{V, k}$ given by the reflexive, transitive extension of $\prec$.
\enddefinition

\proclaim{Lemma \secba.2}
In the notation of Definition \secba.1, $\leq_{e, k}$ is a partial order
on $e_{V, k}$.
\endproclaim

\demo{Proof}
It suffices to show that $\leq_{e, k}$ is antisymmetric.  Suppose for a 
contradiction
that this is not the case; this implies that there is a sequence $$
e_1 \prec e_2 \prec \cdots \prec e_k \prec e_{k+1} = e_1
$$ with $k \geq 2$ and $e_i \ne e_{i+1}$ for all $i$.
Let $d_1, \ldots, d_{k+1}$ be the unique elements
of $d_{V, k}$ for which $(e_i, d_i) \in V_k$.  It then follows that $$
e_1, d_1, e_2, d_2, \ldots, e_k, d_k, e_1
$$ is a nontrivial closed $V$-path, which is the required contradiction.
\qed\enddemo

\proclaim{Proposition \secba.3}
Maintain the notation of Definition \secba.1, and define $$
N = |d_{V, k}| = |e_{V, k}|
.$$  Denote the elements of $d_{V, k}$ 
by $d_1, \ldots, d_N$ in an arbitrary (but fixed) order, and denote by
$e_i$ the element of $e_{V, k}$ paired with $d_i$.
Let $\leq_{e, k}$ be the
partial order on $e_{V, k}$ defined in Lemma \secba.2, and let $\leq_{d, k}$
be the order on $d_{V, k}$ induced by the matching $V$.
\item{\rm (i)}{The matrix
of the linear transformation $\pd_{V, k}$ relative to $d_{V, k}$ and
$e_{V, k}$ is triangular with respect to $\leq_{d, k}$ and $\leq_{e, k}$
with all the diagonal entries equal to $\pm 1$.
In particular, $\pd_{V, k}$ is an isomorphism of $R$-modules.}
\item{\rm (ii)}{Suppose that there exists a $c$ with $1 \leq c \leq N$ 
such that whenever we have $i \leq c < j$, $e_j$ is not a face of $d_i$.
Then if $d_p$ appears with nonzero coefficient in some $d \in D_{V, k}$ for
some $p > c$, then $e_q$ appears with nonzero coefficient in
$\pd(d)$ for some $q > c$.}
\endproclaim

\demo{Proof}
Let $(e, d) \in V_k$.  It follows from the definitions of the partial
orders that $$
\pd_{V, k}(d) = \sum_{{\be \in e_{V, k}} \atop {\be \leq e}} \lambda_\be \be
.$$  It follows from Proposition \secaz.2 that 
$\lambda_\be \in \{-1, 0, +1\}$ for all $\be$, and also that $\lambda_e \ne 0$.
This completes the proof of the first assertion of (i), and 
the second assertion of (i) is immediate from the first.

We now turn to (ii); write $$
d = \sum_{i = 1}^N \lambda_i d_i
.$$  Let $I$ be the set of all $p$ with $c < p \leq N$ satisfying
the hypotheses of (ii), and let $l$ be a $\leq_{d, k}$-maximal element of $I$.
It follows from the definitions that $e_l$ appears with nonzero coefficient
in $\pd(d_l)$.  If $1 \leq i \leq c$, then 
the definition of $c$ shows that $e_l$ appears with zero
coefficient in $\pd(d_i)$, whereas if $c < i \leq N$, then $e_l$ can only
appear with nonzero coefficient in $\pd(d_i)$ if $e_l \leq_{e, k} e_i$, by
the definition of $\pd(d_i)$ and Proposition \secaz.2.  The maximality
hypothesis on $l$ then shows that the only term in the expression $$
\pd(d) = \sum_{i = 1}^N \lambda_i \pd(d_i)
$$ that contributes a coefficient of $e_l$ is the term $\pd(d_l)$.  Setting
$q = l$ completes the proof.
\qed\enddemo

\proclaim{Lemma \secba.4}
Maintain the notation of Definition \secba.1, and suppose that there
are no unpaired $k$-cells.
Suppose also that $e = \sum_{\a \in P^k} \lambda_\a e_\a$ is a 
$k$-cycle; that is, $\pd(e) = 0$.  
\item{\rm (i)}{If $e \ne 0$, then there exists $\a \in P^k$ with $\lambda_\a
\ne 0$ such that $e_\a$ is an upward matching face.}
\item{\rm (ii)}{If
$e = \sum_{\a \in P^k} \lambda_\a e_\a$ and 
$e' = \sum_{\a \in P^k} \mu_\a e_\a$ are two $k$-cycles with the
property that $\lambda_\a = \mu_\a$ whenever $e_\a$ is an upward matching
face, then $e = e'$.}
\endproclaim

\demo{Proof}
Suppose that $e \ne 0$, but that $\lambda_\a = 0$ for
every upward matching face $e_\a \in P^k$.  It follows that $$
e = \sum_{\a \in d_{V, k-1}} \lambda_\a e_\a
.$$  Since $e \ne 0$, Proposition \secba.3 (i) shows 
that $\pd_{V, k-1}(e) \ne 0$.
The definitions of $\pd$ and $\pd_{V, k-1}$ then show 
that $\pd(e) \ne 0$.  This is a contradiction, and (i) follows.  Part (ii)
follows from (i) by considering the cycle $e - e'$.
\qed\enddemo

\proclaim{Lemma \secba.5}
Maintain the notation of Definition \secba.1, and suppose that there
are no unpaired $k$-cells.  Then the set $$
\{ \pd(d) : d \in d_{V, k} \}
$$ is an irredundantly described free $R$-basis for the $k$-cycles over
$R$, $\ker(\pd_k)$.
\endproclaim

\demo{Proof}
Let $e$ be a $k$-cycle.  Since there are no unpaired $k$-cells, we may write $$
e =  \sum_{\a \in d_{V, k-1}} \lambda_\a e_\a + 
\sum_{\be \in e_{V, k}} \mu_\be e_\be
.$$  By Proposition \secba.3 (i), there exists a unique element $$
f = \sum_{\g \in d_{V, k}} \nu_\g e_\g
$$ such that $\pd_{V, k}(f) = 
\sum_{\be} \mu_\be e_\be$.  It follows that for each
$\be \in e_{V, k}$, the coefficient of $e_\be$ in the $k$-cycle $\pd(f)$ is
$\mu_\be$.  Applying Lemma \secba.4 (ii) to the cycles $e$ and $\pd(f)$
then shows that $\pd(f) = e$, and it follows that the set in the statement
is an $R$-spanning set for the cycles.  The freeness assertion follows by
another application of Proposition \secba.3 (i).
\qed\enddemo

\proclaim{Lemma \secba.6}
Let $B$ be a free abelian group on $n$ generators, and let $A$ be
an abelian group generated by $\{a_1, a_2, \ldots, a_n\}$.  If there is
a surjective group homomorphism $\phi : A \ra B$, then $\phi$ is an
isomorphism and $\{a_1, a_2, \ldots, a_n\}$ is a free basis for $A$.
\endproclaim

\demo{Proof}
Since $A$ is generated by $n$ elements, it is naturally 
a homomorphic image $\psi(X)$ of a 
free abelian group $X$ on $n$ generators; it follows that $B = \phi(\psi(X))$
is also a quotient of $X$.

It follows (for example by using the fact that $\kyu$ is a flat $\zed$-module)
that if $$
0 \ra C \ra X \ra B \ra 0
$$ is a short exact sequence of abelian groups, then $\rank(X) = \rank(B)
+ \rank(C)$.  The hypotheses force $\rank(C) = 0$ in this case, but since every
subgroup of $X$ is torsion-free, we must have $C = 0$ and 
$\phi \circ \psi$ is an isomorphism.  It follows that $\psi$ is injective
and is an isomorphism, which completes the proof.
\qed\enddemo

\proclaim{Theorem \secba.7}
Let $X$ be a finite regular CW complex with a $-1$-dimensional cell, and
suppose that $V$ is an acyclic matching on $X$.
Let $Y$ be a CW subcomplex of $X$ and let $V_Y$ be the acyclic partial
matching on $Y$ obtained by discarding all pairs of the matching $V$
that do not entirely lie within $Y$.  Let $K_Y$ be the set of unpaired
cells in $V_Y$, and let $K_X$ be the set of cells of $X \backslash Y$ that
were paired with the elements of $K_Y$ in the original matching $V$.
Suppose that (a) $V$ has no unpaired cells, and that (b) the topological
boundary of each cell of $K_X$ lies entirely within $Y$.  Then
\item{\rm (i)}{the image in $H_k(Y; R)$ of the set $$
{\Cal B}_{Y, k} = \{\pd(k) : k \in K_X \cap P^{k+1} \}
$$ is an $R$-spanning set for $H_k(Y; R)$, and}
\item{\rm (ii)}{if $H_k(Y; \zed)$ is free over $\zed$ of rank $|K_X|$, 
then the image of ${\Cal B}_{Y, k}$ is a free $\zed$-basis for $H_k(Y; \zed)$.}
\endproclaim

\demo{Proof}
To prove (i), we need to show that the map $\pd$ induces a surjective
map from ${\Cal B}_{Y, k}$ to $H_k(Y; R)$.  Let $y$ be a $k$-cycle of $Y$;
we may regard $y$ as a $k$-cycle of $X$ by extension.  Number the elements
of $e_{V, k}$ as $e_1, e_2, \ldots, e_N$ in such a way that $$
e_{V, k} \cap Y = \{e_1, e_2, \ldots, e_c\}
$$ for some $1 \leq c \leq N$.
Since $Y$ is a subcomplex of $X$, we may
choose the numbering so that $$
d_{V, k} \cap Y = \{d_{b+1}, d_{b+2}, \ldots, d_c\}
$$ for some $0 \leq b \leq c$; it follows that $$
K_X = \{d_1, d_2, \ldots, d_b\}
.$$

Hypothesis (b) shows that if $e_j$ is a
face of $d_i$ for $1 \leq i \leq b$, then we must have $j \leq c$.  The
same is true if we have $b < i \leq c$, because $Y$ is a subcomplex
of $X$.  It follows that if $i \leq c < j$, then $e_j$ is not a face
of $d_i$.

Lemma \secba.5 and hypothesis (a) show that
there exists $x \in D_{V, k}$ such that
$\pd(x) = y$; let us write $$
x = \sum_{i = 1}^N \mu_i d_i
.$$   The previous paragraph and Proposition \secba.3 (ii) show that
we must have $\mu_i = 0$ whenever $i > c$; that is, we have $$
x = \sum_{i = 1}^c \mu_i d_i
.$$  If we define $$
x' = \sum_{i = 1}^b \mu_i d_i
,$$ it follows by hypothesis (b) that $\pd(x')$ is a cycle in $Y$.
If $b < i \leq c$ then $d_i$ lies in $Y$; this means that $\pd(d_i)$ is a 
boundary in $Y$ and that $\pd(x)$ and $\pd(x')$ correspond to the same
homology class in $H_k(Y; R)$.  This proves part (i).

In the special case $R = \zed$, we note that the $\zed$-spanning set given 
in (i) has cardinality $|K_X|$.  Part (ii) then follows from Lemma 
\secba.6.
\qed\enddemo

\remark{Remark \secba.8}
The hypotheses of Theorem \secba.7 together with Theorem \secaz.4 (ii)
show that $X$ must be contractible.
\endremark

\head \secc. The half cube \endhead

An $n$-dimensional {\it (Euclidean) polytope} 
$\Pi_n$ is a closed, bounded, convex subset of $\real^n$ obtained by 
intersecting finitely many closed half-spaces associated to hyperplanes.  
We will assume that the set of hyperplanes is taken to be minimal.
The part of $\Pi_n$ that lies in one of
the hyperplanes is called a {\it facet}, and each facet is an 
$(n-1)$-dimensional polytope.  
A polytope is homeomorphic to an $n$-ball (which follows, for example, from
\cite{{\bf 13}, Lemma 1.1}), and the boundary of the polytope,
which is equal to the union of its facets, is identified with the 
$(n-1)$-sphere by this homeomorphism.

Iterating this construction gives
rise to a set of $k$-dimensional polytopes $\Pi_k$ (called {\it $k$-faces}) 
for each $0 \leq k \leq n$.  The elements of $\Pi_0$ are called {\it vertices}
and the elements of $\Pi_1$ are called {\it edges}.  It is not hard to show
that a polytope is the convex hull of its set of vertices, and that the 
boundary of a polytope is precisely the union of its $k$-faces for 
$0 \leq k < n$.  What is less obvious, but still true 
\cite{{\bf 17}, Theorem 1.1}, is that the convex hull of an arbitrary 
finite subset
of ${\Bbb R}^n$ is a polytope in the above sense.  It follows that a polytope
is determined by its vertex set, and we write $\Pi(V)$ for the polytope 
whose vertex set is $V$.  Recalling the vertex sets $\oldpsi_n$, $\ven$ and
$\von$ from the Introduction, we see that $\Pi(\oldpsi_n)$ is an 
$n$-dimensional hypercube, and the half cube $\hcn$ is (by definition) 
$\Pi(\ven)$.

The {\it dimension} of a face is the dimension of its affine hull.  
An {\it automorphism} of a polytope is an isometry of its affine hull that
fixes the polytope setwise.
The {\it interior} of a face refers to its interior with respect to the induced
topology on its affine hull.  

The Coxeter group $W(D_n)$ is a subgroup of the group of geometric 
automorphisms of the half cube $\hcn$, and is the full automorphism group
if $n > 4$.  It is generated by a set of $n$ involutions $\{s_1, s_{1'}, s_2,
s_3, \ldots, s_{n-1}\}$ which act on the set $\oldpsi_n$: the
$(n-1)$ generators $s_i$ act by permuting the coordinates by the transposition
$(i, i+1)$, and the generator $s_{1'}$ acts by the transposition $(1, 2)$
followed by a sign change on the first and second coordinates.  The
group $W(D_n)$ has order $2^{n-1}n!$, and acts on $\oldpsi_n$ by the 
subgroup of signed permutations that effect an even number of sign changes.
This induces an action on $W(D_n)$ on $\real^n$ by orthogonal transformations
fixing the half cube setwise.

\definition{Definition \secc.1}
Let $n \geq 4$ be an integer, and let $\bn = \{1, 2, \ldots, n\}$.

Let $\oldbv' \in \von$ and $S \subseteq \bn$.  We define the subset 
$K(\oldbv', S)$
of $\ven$ by the condition that $\oldbv \in K(\oldbv', S)$ if and only if 
there exists $i \in S$ such that 
$\oldbv$ and $\oldbv'$ differ only in the $i$-th coordinate.

Let $\oldbv \in \ven$ and let $S \subseteq \bn$.
We define the subset $L(\oldbv, S)$ of $\ven$ by the 
condition that $\oldbv' \in L(\oldbv, S)$ if and only if for all $i \not\in S$,
$\oldbv$ and $\oldbv'$ agree in the $i$-th coordinate.  The set $S$ is 
characterized as the set of coordinates
at which not all points of $L(\oldbv, S)$ agree.

We call the set $S$ in a face of the form $\Pi(K(\oldbv', S))$ or
$\Pi(L(\oldbv, S))$ the {\it mask} of the face.
\enddefinition

The $k$-faces of the half cube were classified in \cite{{\bf 10}}.

\proclaim{Theorem \secc.2 \cite{{\bf 10}}}
The $k$-faces of $\hcn$ for $k \leq n$ are as follows:
\item{\rm (i)}{$2^{n-1}$ $0$-faces (vertices) given by the elements of
$\ven$;}
\item{\rm (ii)}{$2^{n-2} {n \choose 2}$ $1$-faces $\Pi(K(\oldbv', S))$,
where $\oldbv' \in \von$ and $|S| = 2$;}
\item{\rm (iii)}{$2^{n-1} {n \choose 3}$ simplex shaped $2$-faces 
$\Pi(K(\oldbv', S))$, where $\oldbv' \in \von$ and $|S| = 3$;}
\item{\rm (iv)}{$2^{n-1} {n \choose {k+1}}$ simplex shaped $k$-faces 
$\Pi(K(\oldbv', S))$, where $\oldbv' \in \von$ and $|S| = k+1$ for 
$3 \leq k < n$;}
\item{\rm (v)}{$2^{n-k} {n \choose k}$ half cube shaped $k$-faces 
$\Pi(L(\oldbv, S))$, where $\oldbv \in \ven$ and $|S| = k$ for 
$3 \leq k \leq n$.}

Furthermore, two faces are conjugate under the action of $W(D_n)$ if and
only if they have the same dimension and the same shape.
\endproclaim

\demo{Proof}
The classification of the $k$-faces is given in \cite{{\bf 10}, Theorem 2.3.6},
and the classification of the orbits under the action of $W(D_n)$ is given
in \cite{{\bf 10}, Theorem 4.2.3 (ii)}.
\qed\enddemo

The unique $n$-face in (v) above corresponds to the interior of the polytope.
The $k$-faces assemble naturally into a regular CW complex, $C_n$.  

\remark{Remark \secc.3}
The proof of Theorem \secc.2 given in \cite{{\bf 10}} is not optimal.  A
shorter proof of this result can be obtained by using Casselman's theorem
\cite{{\bf 5}, Theorem 3.1}, which Casselman attributes to Satake and Borel--Tits.
The latter result gives an explicit set of $W(D_n)$-orbit representatives of 
the $k$-faces of the half cube for each $k$.
\endremark

In order to describe a complete acyclic matching on the faces of the half
cube, it will be helpful to parametrize the faces in terms of certain
sequences.  

\definition{Definition \secc.4}
We denote a coordinate of $+1$ by the digit $0$, and a coordinate
of $-1$ by the digit $1$.  A face of type $\Pi(K(\oldbv', S))$ is denoted by
replacing the digits in $v'$
corresponding to coordinates in $S$ by underlined symbols.
A face of type $\Pi(L(\oldbv, S))$ is denoted by replacing the digits in $v$
corresponding to $S$ by asterisks.  This notation is ambiguous for the
$1$-dimensional faces; we consider them to be faces of type 
$\Pi(K(\oldbv', S))$ in which the symbol associated to the rightmost coordinate
in $S$ is a $\boz$.
\enddefinition

\example{Example \secc.5}
\item{(i)}{The vertex $(-1, -1, -1, +1, -1, +1, +1)$ corresponds to the
sequence $1110100$.}
\item{(ii)}{The simplex shaped face $\Pi(K(\oldbv', S))$ with $$
\oldbv' = (+1, -1, -1, +1, -1, +1, +1)
$$ and $S = \{1, 3, 6, 7\}$ is denoted by the sequence $\boz1\boo01\boz\boz$.
By toggling each coordinate in $S$ in turn, we find the set of vertices of
this face; these vertices correspond to the sequences 
$1110100$,
$0100100$, 
$0110110$, 
$0110101$.}
\item{(iii)}{The half cube shaped face $010{*}{*}1{*}010$ is the convex 
hull of the $2^3$ points
obtained by filling in the asterisks with $0$s and $1$s in such a way that
the total number of $1$s is even, \idest
$0100011010$,  
$0100110010$,  
$0101010010$,  
$0101111010$.  
This face is equal to $\Pi(L(\oldbv, S))$, where $\oldbv$ is any of the
four points corresponding to the sequences listed, and $S = \{4, 5, 7\}$.}
\item{(iv)}{The convex hull of the pair of vertices $1110100$ and 
$0100100$ is a $1$-dimensional face.  This could potentially be denoted
either by $\boo1\boz0100$ or by $\boz1\boo0100$, but only the first of
these is correct according to Definition \secc.4.

Similarly, the convex hull of the pair of vertices $0110110$ and $0110101$
is denoted by $01101\boz\boz$, rather than $01101\boo\boo$.}
\endexample

We may now describe an explicit matching on the faces of the half cube, 
together with the empty face, $\emptyset$.  Let $F$ be one of these faces,
let $S$ be its mask, and let $\bx$ be the sequence associated to $F$ by
Definition \secc.4.  We denote the face matched with $F$ by $F'$, and the
sequence associated with $F'$ by $\boy$.

\item{(1)}{If $F$ is a half cube shaped face with $\dim(F) \geq 3$, 
and $\bx$ contains a $1$ to the right of $S$, then $\boy$ is obtained by
replacing the rightmost $1$ in $\bx$ with a ${*}$.}

\item{(2)}{If $F$ is a half cube shaped face with $\dim(F) \geq 4$, and there
is no $1$ in $\bx$ to the right of $S$, then $\boy$ is obtained by replacing
the rightmost ${*}$ in $S$ with a $1$.}

\item{(3)}{If $F$ is a simplex shaped face with $\dim(F) \geq 2$, and
the rightmost $1$ in $\bx$ is not underlined, then $\boy$ is obtained by
underlining the rightmost $1$ in $\bx$.}

\item{(4)}{If $F$ is a simplex shaped face with $\dim(F) \geq 3$, and
the rightmost $1$ in $\bx$ is underlined, then $\boy$ is obtained by replacing 
the rightmost $\boo$ with a $1$.}

\item{(5)}{If $F$ is a triangle (a simplex shaped face with $\dim(F) = 2$),
and the rightmost $1$ in $\bx$ is underlined, 
and the entries in $S$ (reading left to right) are $\boz\boo\boo$ or 
$\boo\boo\boo$, then $\boy$ is obtained by replacing these three entries in
by ${*}$.}

\item{(6)}{If $F$ is a half cube shaped face with $\dim(F) = 3$, 
and there is no $1$ to the right of $S$ in $\bx$, then $\boy$ is obtained
from $\bx$ by replacing the rightmost two ${*}$ in $\bx$ by $\boo$, and
replacing the leftmost ${*}$ in $\bx$ by $\boz$ or by $\boo$, in such a
way that the total number of $1$s in $\boy$ is odd.}

\item{(7)}{If $F$ is an edge and the rightmost $1$ in $F$
is not underlined, 
then $\boy$ is obtained from $\bx$ by underlining the rightmost $1$.}

\item{(8)}{If $F$ is a triangle and the rightmost $1$ in $\bx$ is underlined,
and it is not the case that the rightmost two entries in $S$ are equal to 
$\boo\boo$, then $\boy$ is obtained from $\bx$ by replacing the rightmost
$\boo$ with a $1$.}

\item{(9)}{If $F$ is a vertex and $\bx$ contains at least two $1$s, then 
$\boy$ is obtained from $\bx$ by replacing the rightmost $1$ in $\bx$ by 
$\boz$ and the second rightmost $1$ in $\bx$ by $\boo$.}

\item{(10)}{If $F$ is an edge and the rightmost $1$ in $F$
is underlined, then $\boy$ is obtained from $\bx$ by replacing both underlined
entries by non-underlined $1$s.}

\item{(11)}{The empty face $\emptyset$ is matched with the vertex
$00\cdots0$.}

\example{Example \secc.6}

\item{(i)}{The faces $0{*}{*}1{*}10$ and $0{*}{*}1{*}{*}0$ 
are matched by rules (1) and (2).}
\item{(ii)}{The faces $0\boz1\boo10\boz$ and $0\boz1\boo\boo0\boz$ are
matched by rules (3) and (4).}
\item{(iii)}{The faces $0\boo1\boo10\boo$ and $0{*}1{*}10{*}$ are matched
by rules (5) and (6).}
\item{(iv)}{The faces $01\boo01\boz0$ and $01\boo0\boo\boz0$ are matched
by rules (7) and (8).}
\item{(v)}{The faces $1110010$ and $11\boo00\boz0$ are matched by rules
(9) and (10).}
\endexample

\proclaim{Lemma \secc.7}
Every face (including the empty face) of the half cube $\hcn$ is matched
to another face by one, and only one, of rules (1)--(11) above.  Furthermore,
the smaller face in each pair of matched faces is a codimension $1$ face
of the larger of the pair.
\endproclaim

\demo{Proof}
This is a case analysis based on the classification of Theorem \secc.2.
Let $F$ be a (possibly empty) face of $\hcn$, and let $\bx$ and $S$ be the 
corresponding sequence and mask, respectively.

The sequences corresponding to vertices all contain an even number of $1$s.
If this number is zero, then the vertex is matched to $\emptyset$ by (11);
otherwise, the vertex is matched to an edge by (9).  Notice that if rule 
(9) is applied, the resulting sequence satisfies the conditions of
Definition \secc.4.

The sequences corresponding to edges all contain an odd total number of $1$s,
and have rightmost underlined entry equal to $\boz$ by Definition \secc.4.
In particular, there must be at least one $1$ (underlined or otherwise) in
the sequence.  If the rightmost $1$ is not underlined, then $F$ is matched to
a triangle by (7); otherwise, $F$ is matched to a vertex by (10).  Notice
that rule (10) in this case will produce a vertex with an even number of $1$s, 
as required.

If $F$ is a triangle, then $\bx$ contains an odd (and thus nonzero) number of
$1$s, some of which may be underlined.  If the rightmost $1$ is not underlined,
then $F$ is matched to a $3$-simplex by (3).  If the rightmost $1$ is
underlined and the two rightmost underlined symbols are both $1$s, then
$F$ is matched to a $3$-half cube by (5); otherwise, $F$ is matched to an
edge by (8).  If rule (8) is applied, the rightmost underlined symbol in the
resulting edge cannot be a $1$, or rule (5) would have applied instead; this
satisfies the requirements of Definition \secc.4.

If $F$ is a simplex of dimension at least $3$, then $\bx$ contains an odd
number of $1$s; in particular, there is at least one occurrence of $1$.
If the rightmost such occurrence is not underlined, then $F$ is matched to
a simplex of dimension one larger by (3); if the rightmost such occurrence 
is not underlined, then $F$ is matched to a simplex of dimension one less
by (4).

If $F$ is a half cube, and there is a $1$ in $\bx$ to the right of $S$,
then $F$ is matched to a higher dimensional half cube by (1).  Suppose
there is no such $1$.  If $F$ has dimension at least $4$ (respectively,
dimension equal to $3$), then $F$ is matched to a half cube of dimension one 
lower by (2) (respectively, (6)).
\qed\enddemo

\proclaim{Lemma \secc.8}
Let $F_1$ and $F_2$ be faces of $\hcn$.  Rule (1) (respectively,
(3), (5), (7), (9)) sends face $F_1$ to the face $F_2$ if and only if
rule (2) (respectively, (4), (6), (8), (10)) sends $F_2$ to $F_1$.
\endproclaim

\demo{Proof}
It is immediate from the definitions that rules (1) and (2) are inverses of
each other, restricted to the appropriate domain and codomain.  The same is
true for rules (3) and (4).

Observe that if $F$ is a $3$-dimensional half cube shaped face, then $F$ 
contains four triangular faces.  These are obtained by replacing the $*$ in
the mask of $S$ by occurrences of $\boo$ or $\boz$ in such a way that the
total number of $1$s and $\boo$s in the resulting sequence is odd.
Precisely one of these four triangular faces has a mask of the form
$\boz\boo\boo$ or $\boo\boo\boo$.  These observations imply that rules (5)
and (6) are inverses of each other.

Note that if $F$ is an edge, then Definition \secc.4 requires the rightmost
underlined symbol in $F$ to be a $\boz$.  If rule (7) is applicable to $F$,
and $S$ is the mask of the resulting triangle, then the rightmost two entries
in $S$ will be $\boz\boo$.  On the other hand, if rule (8) is applicable
to a triangle $F'$ with mask $S'$, then the rightmost two entries in $S$
will be $\boz\boo$, and the rightmost underlined entry of the resulting edge
will be $\boz$.  These observations show that rules (7) and (8) are inverses
of each other.

Observe that if $F$ is an edge and the rightmost $1$ in $F$ is underlined,
then the restrictions of Definition \secc.4 mean that this rightmost $1$ is
the leftmost of the two underlined symbols, and furthermore, that the 
rightmost underlined symbol is a $\boz$.  Since $F$ contains an odd number
of occurrences of $1$ or $\boo$, replacing both these entries with
occurrences of $1$ as in rule (10) will produce an even total number of $1$s.
Conversely, any vertex not equal to $00\cdots 0$ contains an even number
of $1$s; in particular, it contains at least two occurrences.  These 
observations show that rules (9) and (10) are inverses to each other, 
which completes the proof.
\qed\enddemo

The following result is an immediate consequence of lemmas 
\secc.7 and \secc.8.

\proclaim{Proposition \secc.9}
Rules (1)--(11) define a complete matching on the set of faces of the half
cube $\hcn$, including the empty face.
\qed\endproclaim

\head \secd. Proof that the matching is acyclic \endhead

In \S\secd, we will show that the complete matching of Proposition \secc.9
is acyclic in the sense of \S\secaz.  In order to do this, it is convenient
to introduce a certain statistic on the faces of types (i)--(iv) in the
classification of Theorem \secc.2; we will call such faces {\it faces of
type $K$}.

\definition{Definition \secd.1}
Let $F$ be a face of type $K$, and let $\bs$ be the sequence associated to
$F$ by Definition \secc.4.  Let $S'
\subseteq \bn$ be the (possibly empty) set of indices at which $\bs$ has 
occurrences of $1$ or $\boo$.
We define the {\it total} of $F$ to be $$
t(F) = t(\bs) = \sum_{i \in S'} i.
$$  We define the sequence $u(\bs) = u(F)$ from $\bs$ by replacing all
occurrences of $\boz$ (respectively, $\boo$) by $0$ (respectively, $1$).
\enddefinition

It is immediate $u(F) = u(F')$ implies that $t(F) = t(F')$.

\example{Example \secd.2}
If $F$ is the face with sequence $0\boo11\boz\boo01$, then
we have $t(F) = 2 + 3 + 4 + 6 + 8 = 23$ and $u(F) = 01110101$.
\endexample

\remark{Remark \secd.3}
\item{(i)}{In Definition \secc.4, the sequence chosen to represent a
given edge is the one with the lower of the two possible totals.}
\item{(ii)}{In the context of rule (6) of the matching, there are four
triangular faces of the $3$-dimensional half cube; the one paired with
the half cube is the one with the highest total.}
\endremark

The following result is a immediate from the classification of Theorem
\secc.2; it will often be used in the sequel.

\proclaim{Lemma \secd.4}
\item{\rm (i)}{If $\Pi(K(\oldbv', S))$ is a simplex shaped face of $\hcn$, 
then every face of $\Pi(K(\oldbv',S))$ can be expressed in the form 
$\Pi(K(\oldbv',S'))$ for the same $\oldbv'$, where $S' \subset S$ and
$S \backslash S'$ is a singleton.}
\item{\rm (ii)}{If $\Pi(L(\oldbv', S))$ is a half cube shaped face of $\hcn$, 
then every half cube shaped face of $\Pi(L(\oldbv',S))$ can be expressed 
in the form $\Pi(L(\oldbv',S'))$ for the same $\oldbv'$, where 
$S' \subset S$ and $S \backslash S'$ is a singleton.}
\qed\endproclaim

\remark{Remark \secd.5}
Some care must be taken in using Lemma \secd.4 (i) for faces of the form
$\Pi(K(\oldbv',T))$ when $|T| = 2$, because in this case, there are
two possible choices for $\oldbv'$.
\endremark

\proclaim{Lemma \secd.6}
If $$
\a_0, \be_0, \a_1, \be_1, \a_2, \ldots, \be_r, \a_{r+1} = \a_0
$$ is a nontrivial 
closed $V$-path of faces of $\hcn$ in which the faces $\a_i$ have
dimension $2$, then all the faces $\a_i$ and $\be_i$ are of the form
$\Pi(K(\oldbv', S))$ for the same $\oldbv'$.  
In particular, none of the $\be_i$ is half cube shaped, and the 
sequences $u(\a_i)$ and $u(\be_i)$ all coincide.
\endproclaim

\demo{Proof}
If a face $\be_i$ is a $3$-dimensional half cube,
it follows from Remark \secd.3 (ii) that $t(\a_{i+1}) < t(\a_i)$.  In
contrast, if $\be_i$ is a $3$-dimensional simplex, Lemma \secd.4 (i)
shows that $u(\a_i) = u(\be_i) = u(\a_{i+1})$, which in turn implies that 
$t(\a_{i+1}) = t(\a_i)$.  The conclusions (i)
and (ii) now follow from the requirement that $t(\a_{r+1}) = t(\a_0)$.
\qed\enddemo

\proclaim{Lemma \secd.7}
If $$
\a_0, \be_0, \a_1, \be_1, \a_2, \ldots, \be_r, \a_{r+1} = \a_0
$$ is a nontrivial closed $V$-path of faces of $\hcn$ in which the faces 
$\a_i$ have
dimension $1$, then all the sequences $u(\a_i)$ and $u(\be_i)$ all coincide.
\endproclaim

\demo{Proof}
It follows from Remark \secd.5 and Remark \secd.3 (i) that we have
$t(\a_{i+1}) \leq t(\a_i)$ for $0 \leq i \leq r$, with equality holding
if and only if $u(\a_i) = u(\be_i) = u(\a_{i+1})$.  The requirement that
$\a_{r+1} = \a_0$ forces equality to hold at every step, and the conclusion
follows from this.
\qed\enddemo

\proclaim{Theorem \secd.8}
The matching described in \S\secc\  is a complete acyclic matching on the
faces of $\hcn$ (together with the empty face).
\endproclaim

\demo{Proof}
By Proposition \secc.9, it is enough to show that the matching is
acyclic.
By Theorem \secaz.4 (i), this reduces to showing that there are no nontrivial
closed $V$-paths.  Suppose for a contradiction that $$
\a_0, \be_0, \a_1, \be_1, \a_2, \ldots, \be_r, \a_{r+1} = \a_0
$$ is such a path.  We will proceed by a case analysis based on
$\dim(\a_0)$.

The fact that there is a unique face of dimension $-1$ rules out the
possibility of $\dim(\a_0) = -1$.

Suppose that $\dim(\a_0) = 0$.  Each edge $\be_i$ has exactly two vertices
contained in it, and they both appear in the path.  It follows from 
Remark \secd.3 (i) that for all $0 \leq i \leq r$, $t(\a_i) < t(\a_{i+1})$.
This is incompatible with the condition $\a_{r+1} = \a_0$, which completes
the proof in this case.

Suppose that $\dim(\a_0) \geq 3$ and that $\a_0$ is simplex shaped.  It
follows by rule (3) of the matching and Lemma \secd.4 (i) that all the other
faces in the path are simplex shaped, with all the matched pairs being
matched by rules (3) and (4).  In particular, $\be_0$ is obtained by 
underlining the rightmost $1$ in the sequence for $\a_0$, and $\a_1$ is
obtained from $\be_0$ by removing the underline from one of the other symbols.
This means that the rightmost $1$ in the sequence of $\a_1$ is still 
underlined, and $\a_1$ is not a candidate for input to rule (3).  This
is a contradiction.

If $\a_0$ is a triangle, Lemma \secd.6 shows that none of the $\be_i$ is
half cube shaped.  In particular, rules (5) and (6) do not play a role in
the path, and we can apply the argument of the above paragraph to obtain
a contradiction.

If $\a_0$ is an edge, then all the matched pairs in the path are matched by
rules (7) and (8).  By Lemma \secd.7, the sequences $u(\a_i)$ and $u(\be_i)$
all coincide.  We can then copy the argument used above to deal with the
case where $\a_i$ is a simplex to obtain a contradiction.

It remains to deal with the case where $\a_0$ is a half cube shaped face
with $\dim(\a_0) \geq 3$.  If at least one of the $\a_i$ is a simplex shaped
face, then we may rotate the closed path so that $\a_i$ plays the role of
$\a_0$.  This has already been dealt with above, so we may assume that
all the $\a_i$ are half cube shaped, and that all faces in the path are 
matched by rules (1) and (2).

It follows from rule (1) that $\be_0$ is obtained by 
replacing the rightmost $1$ by a $*$ in the sequence for $\a_0$, and $\a_1$ is
obtained from $\be_0$ via Lemma \secd.4 (ii) 
by replacing one of the other occurrences of $*$ by
$0$ or $1$.  This means that $\a_1$ has no $1$ to the 
right of the rightmost $*$,
and is not a candidate for input to rule (1).  This is a contradiction and
completes the proof.
\qed\enddemo

\head \sece. Homology bases for polytopal subcomplexes \endhead

In this section, we combine Theorem \secd.8 with Theorem \secba.7 to 
obtain an explicit homology basis for $C_{n, k}$.

\proclaim{Lemma \sece.1}
Let $n \geq 4$ and let $3 \leq k < n$.  Let $X$ be the CW complex 
corresponding to the faces of $\hcn$, including the empty face, let $V$
be the complete acyclic matching on $X$ given in Theorem \secd.8, and let $Y$
be the subcomplex of $X$ corresponding to $C_{n, k}$.  Then:
\item{\rm (i)}{$X$ and $Y$ satisfy the hypotheses of Theorem \secba.7;}
\item{\rm (ii)}{the unmatched faces of $Y$ are the $(k-1)$-dimensional faces
that are inputs to rules (1) or (5) of the matching; these are paired
with the $k$-dimensional half cube shaped faces of $X$ that are inputs to
rules (2) or (6).}
\endproclaim

\demo{Proof}
We need to identify the faces of $Y$ that are paired by $V$ with faces
in $X \backslash Y$.  An inspection of the matching rules in \S\secc\ shows
that these faces are the faces of $Y$ of dimension $k-1$ satisfying the
input conditions of rule (1) if $k > 3$, or rule (5) if $k = 3$.
The faces of  $X \backslash Y$ that are paired with these faces are 
$k$-dimensional half
cubes that satisfy the input conditions of rule (2) if $k > 3$, or rule
(6) if $k = 3$; this proves
(ii).  

The faces of a $k$-dimensional half cube shaped face are 
$(k-1)$-dimensional, and are 
either simplices or half cubes.  All such faces are contained
in $Y$.  This shows that condition (b) of Theorem \secba.7 holds, and 
condition (a) holds by the completeness of the matching $V$, completing the
proof of (i).
\qed\enddemo

\proclaim{Theorem \sece.2}
Let $n \geq 4$ and $3 < k < n$, and 
let $B$ be the set of $k$-dimensional half cube shaped faces of $\hcn$ 
whose sequences have no $1$ to the right of the rightmost occurrence of $*$.
\item{\rm (i)}{A basis for the $(k-1)$-st homology of $C_{n, k}$ is given 
by the images under the boundary map of the faces in $B$.}
\item{\rm (ii)}{The $(k-1)$-st Betti number of $C_{n, k}$ is given by $$
\sum_{i = 1}^n 2^{i-k} {{i-1} \choose {k-1}}. \eqno{(\secz.2)}
$$}
\endproclaim

\demo{Proof}
The hypotheses of Theorem \secba.7 are satisfied by Lemma \sece.1 (i).  
By Lemma \sece.1 (ii), the set ${\Cal B}_{Y, k}$ in Theorem \secba.7 consists
of the $k$-dimensional half cube shaped faces that are inputs to rules (2)
or (6), and the latter coincides with the set $B$ by the definition of the
matching.

The $(k-1)$-st reduced homology of $C_{n, k}$ is free over $\zed$ by 
\cite{{\bf 10}, Theorem 3.3.2}.  Part (i) now follows from Theorem 
\secba.7 (ii).

For part (ii), notice that the faces of $B$ all have sequences with 
precisely $k$ occurrences of $*$, and furthermore, they all end in
${*}00\cdots0$.  Let $i$ denote the number of symbols including and to the 
left of the rightmost $*$, so that the number of symbols in the 
sequence ${*}00\cdots0$ just mentioned is $n-i+1$.  To form the set of
such sequences for a fixed $i$, the leftmost $i-1$ symbols must contain
$k-1$ occurrences of $*$; the remaining $i-k$ symbols can be independently
chosen to be $0$ or $1$.  (This number will be zero unless $i \geq k$.)  
This gives a total of $2^{i-k} {{i-1} \choose {k-1}}$
choices, and summing over all possible $i$ gives the result.
\qed\enddemo

\remark{Remark \sece.3}
The basis of Theorem \sece.2 can be used for explicit 
computations involving the action of $W(D_n)$ on the integral homology of
$C_{n,k}$.  In this case, the incidence numbers may be computed using the 
combinatorics of Coxeter groups.
\endremark

\head Acknowledgements \endhead

This material is based upon work supported by the National Science Foundation 
under Grant Number 0905768.

Any opinions, findings, and conclusions or recommendations expressed in this 
material are those of the authors and do not necessarily reflect the views of 
the National Science Foundation.

\leftheadtext{} \rightheadtext{}

\end
\Refs\refstyle{A}\widestnumber\key{{\bf 13}} \leftheadtext{References}
\rightheadtext{References} 

\ref\key{{\bf 1}}
\by H. Barcelo and S. Smith
\paper The discrete fundamental group of the order complex of $B_n$
\jour J. Algebraic Combin.
\vol 27 \yr 2008 \pages 399--421
\endref

\ref\key{{\bf 2}}
\by A. Bj\"orner
\paper Random walks, arrangements, cell complexes, greedoids, 
and self-organizing libraries
\inbook Building Bridges
\publ Springer-Verlag
\publaddr Berlin
\pages 165--203
\vol 19
\yr 2008
\endref

\ref\key{{\bf 3}}
\by A. Bj\"orner and V. Welker
\paper The homology of ``k-equal'' manifolds and related partition lattices
\jour Adv. Math.
\vol 110 \yr 1995 \pages 277--313
\endref

\ref\key{{\bf 4}}
\by H. Bruggesser and P. Mani
\paper Shellable decompositions of cells and spheres
\jour Math. Scand.
\vol 29 \yr 1971
\pages 197--205
\endref

\ref\key{{\bf 5}}
\by W. Casselman
\paper Geometric rationality of Satake compactifications
\inbook Algebraic groups and Lie groups
\publ Cambridge University Press
\publaddr Cambridge, UK
\pages 81--103
\yr 1997
\endref

\ref\key{{\bf 6}}
\by M. Chari
\paper On discrete Morse functions and combinatorial decompositions
\jour Discrete Math.
\vol 217 \yr 2000 \pages 101--113
\endref

\ref\key{{\bf 7}}
\by R. Forman
\paper Morse Theory for cell complexes
\jour Adv. Math.
\vol 134 \yr 1998 \pages 90--145 
\endref

\ref\key{{\bf 8}}
\by R. Forman
\paper A user's guide to discrete Morse theory
\jour S\'eminaire Lotharingien de Combinatoire
\vol 48 \yr 2002
\endref

\ref\key{{\bf 9}}
\by R. Geoghegan
\book Topological Methods in Group Theory
\publ Springer
\publaddr New York
\yr 2008
\endref

\ref\key{{\bf 10}}
\by R.M. Green
\paper Homology representations arising from the half cube
\jour Adv. Math.
\vol 222 \yr 2009 \pages 216--239
\endref

\ref\key{{\bf 11}}
\by R.M. Green
\paper Homology representations arising from the half cube, II
\jour J. Combin. Theory Ser. A
\vol 117 \yr 2010 \pages 1037--1048
\endref

\ref\key{{\bf 12}}
\by J.T. Harper
\book Homology representations arising from a hypersimplex
\publ Ph.D. thesis
\publaddr University of Colorado Boulder
\yr 2011
\endref

\ref\key{{\bf 13}}
\by J.R. Munkres
\book Elements of algebraic topology
\publ Addison-Wesley
\publaddr Menlo Park, CA
\yr 1984
\endref

\ref\key{{\bf 14}}
\by J.J. Rotman
\book An introduction to algebraic topology
\publ Springer-Verlag
\publaddr New York
\yr 1988
\endref

\ref\key{{\bf 15}}
\by M. Shattuck and T. Waldhauser
\paper Proofs of some binomial identities using the method of last squares
\miscnote preprint; {\tt arXiv:1107.1063}
\endref

\ref\key{{\bf 16}}
\by N.J.A. Sloane
\paper The on-line encyclopedia of integer sequences
\miscnote available online at 
\newline {\tt http://www.research.att.com/~njas/sequences/}
\yr 2011
\endref

\ref\key{{\bf 17}}
\by G.M. Ziegler
\book Lectures on polytopes
\publ Springer-Verlag
\publaddr New York
\yr 1995
\endref

\endRefs
